\documentclass[12pt]{article}
\usepackage{amssymb,amsmath,bm}
\usepackage{mathrsfs}
\usepackage{hyperref}
\usepackage{graphicx}
\hypersetup{
    colorlinks,
    citecolor=black,
    filecolor=black,
    linkcolor=black,
    urlcolor=black
}
\usepackage{enumerate}
\begin{document}
\newcommand{\beq}{\begin{eqnarray}}
\newcommand{\eeq}{\end{eqnarray}}
\newcommand{\beas}{\begin{eqnarray*}}
\newcommand{\enas}{\end{eqnarray*}}
\newcommand{\bea}{\begin{eqnarray}}
\newcommand{\ena}{\end{eqnarray}}
\newcommand{\bms}{\begin{multline*}}
\newcommand{\ems}{\end{multline*}}
\newcommand{\qmq}[1]{\quad \mbox{#1} \quad}
\newcommand{\qm}[1]{\quad \mbox{#1}}
\newcommand{\nn}{\nonumber}
\newcommand{\bbox}{\hfill $\Box$}
\newcommand{\ignore}[1]{}
\newcommand{\tr}{\mbox{tr}}
\newcommand{\Bvert}{\left\vert\vphantom{\frac{1}{1}}\right.}
\newtheorem{theorem}{Theorem}[section]
\newtheorem{corollary}{Corollary}[section]
\newtheorem{conjecture}{Conjecture}[section]
\newtheorem{proposition}{Proposition}[section]
\newtheorem{remark}{Remark}[section]
\newtheorem{lemma}{Lemma}[section]
\newtheorem{example}{Example}[section]
\newtheorem{definition}{Definition}[section]
\newtheorem{condition}{Condition}[section]
\newcommand{\pf}{\noindent {\bf Proof:} }
\newcommand{\sbull}{\scalebox{0.5}{\textbullet}}
\def\blfootnote{\xdef\@thefnmark{}\@footnotetext}
\newcommand{\lcolor}[1]{\textcolor{magenta}
	{#1}}  
\newcommand{\lcomm}[1]{\marginpar{\tiny\lcolor{#1}}}  

\title{{\bf\Large Non-asymptotic distributional bounds for the Dickman approximation of the running time of the Quickselect algorithm}}

\author{Larry Goldstein}
\date{Department of Mathematics, University of Southern California}

\footnotetext{This work was partially supported by NSA grant
	H98230-15-1-0250.}


\maketitle

\begin{abstract}
Given a non-negative random variable $W$ and $\theta>0$, let the generalized Dickman transformation map the distribution of $W$ to
that of 
\begin{eqnarray*}
	W^*=_d U^{1/\theta}(W+1),
\end{eqnarray*}
where $U \sim {\cal U}[0,1]$, a uniformly distributed variable on the unit interval, independent of $W$, and where $=_d$ denotes equality in distribution. It is well known that $W^*$ and $W$ are equal in distribution if and only if $W$ has the generalized Dickman distribution ${\cal D}_\theta$. We demonstrate that the Wasserstein distance $d_1$ between $W$, a non-negative random variable with finite mean, and $D_\theta$ having distribution ${\cal D}_\theta$ obeys the inequality  
\begin{eqnarray*}
	d_1(W,D_\theta) \le (1+\theta)d_1(W,W^*).
\end{eqnarray*}
The specialization of this bound to the case $\theta=1$ and coupling constructions yield
\begin{eqnarray*}
	d_1(W_{n,1},D_1) \le \frac{8\log (n/2)+10}{n} \quad \mbox{for all $n \ge 1$, where for $m \ge 1$} \quad W_{n,m}=\frac{1}{n}C_{n,m}-1,
\end{eqnarray*}
and $C_{n,m}$ is the number of comparisons made by the Quickselect algorithm to find the $m^{th}$ smallest element of a list of $n$ distinct numbers. A similar bound holds for $W_{n,m}$ for $m \ge 2$, and together recover and quantify the results of \cite{MR1918722} that show distributional convergence of $W_{n,m}$ to the standard Dickman distribution ${\cal D}_1$ in the asymptotic regime $m=o(n)$. By comparison to an exact expression for the expected running time $E[C_{n,m}]$, lower bounds are provided that show the rate is not improvable for $m \not = 2$.
\end{abstract}

\section{Introduction}
For a given non-negative random variable $W$ and $\theta>0$, let the generalized Dickman transformation map the distribution of $W$ to that of
\bea \label{eq:W*transform}
W^*=_d U^{1/\theta}(W+1),
\ena
where $U$ has the uniform distribution ${\cal U}[0,1]$ on the unit interval, and is independent of $W$ and where $=_d$ denotes equality in distribution. It is well known \cite{MR2575452}, \cite{MR2079909} that the generalized Dickman distribution ${\cal D}_\theta$ is the unique fixed point of the  transformation \eqref{eq:W*transform}, that is,
\bea \label{eq:D.is.U(D+1)}
W \sim {\cal D}_\theta \qmq{if and only if} W=_d W^*.
\ena
When \eqref{eq:W*transform} holds we will say that $W^*$ has the ${\cal D}_\theta$-bias distribution of $W$. In what follows, $D_\theta$ will denote a random variable with distribution ${\cal D}_\theta$.
The case $\theta=1$ corresponds to the (standard) Dickman distribution, for which we may drop the subscript $\theta$.

The Dickman function $\rho$ first made its appearance in number theory \cite{dickman1930frequency} when counting the number of integers below a fixed threshold whose prime factors satisfy some given upper bound. Standardizing $\rho$ yields the density of the standard Dickman distribution, the cannonical member of the family ${\cal D}_\theta, \theta>0$ of generalized Dickman distributions, which also arise in the study of component counts of 
logarithmic combinatorial structures such as permutations and partitions \cite{MR2032426}, and more generally for the quasi-logarithmic class considered in \cite{MR2793242}.
See also the recent work \cite{pinsky2016natural}, \cite{azmoodeh2016distances} and \cite{bhattacharjee2017dickman} in this area, that detail some connections to probabilistic number theory. 

Members from the generalized Dickman family have subsequently been noted to arise in a variety of other contexts, in particular for the sum of edge lengths of vertices connected to the origin in minimal directed spanning trees in \cite{MR2079909}, and for weighted sums of independent random variables in \cite{MR3751078}, \cite{azmoodeh2016distances} and \cite{bhattacharjee2017dickman}. Simulation of the Dickman distribution has been considered in \cite{MR2575452}.

Here we study the error incurred when using the standard Dickman distribution to approximate that of the (properly normalized) number of comparisons made by the Quickselect sorting algorithm of Hoare \cite{hoare1961algorithm} for locating the $m^{th}$ smallest element of a list of $n$ distinct numbers. 
One may visualize how Quickselect works in terms of a tree structure. First, a `pivot'  is chosen uniformly from the given list. The list is then divided into those numbers on the list that are strictly smaller, making up the left subtree, and those that are strictly larger, making up the right. If the left subtree is of size $m-1$ then the pivot is the desired $m^{th}$ smallest element, and the procedure terminates. Otherwise, the process continues recursively on the left sub-tree if it is of size $m$ or larger, and else on the right sub-tree.

Letting
\bea \label{eq:intro.def.wnm}
W_{n,m} =\frac{1}{n}C_{n,m}-1,
\ena
where $C_{n,m}$ is the number of comparisons made by Quickselect, the work of \cite{MR1918722} showed that $W_{n,m}$ converges in distribution to the Dickman $D$ when $m=o(n)$. We note that in the case $m=1$ Quickselect simplifies in that at each step of the recursion the procedure either stops or continues on the left subtree. As this case is simpler than for $m \ge 2$ we deal with it separately.

The following two theorems quantify and recover the results of \cite{MR1918722} by providing non-asypmptotic bounds in the Wasserstein distance $d_1$ between $W_{n,m}$ and $D$ that converge to zero in the $m=o(n)$ asymptotic regime. As the $m^{th}$ smallest number of a list of $n$ distinct numbers only exists when $n \ge m$, we need only consider this range of parameters in what follows. 

\begin{theorem} \label{thm:Quickselect}
	Let $C_{n,1}$ be the number of comparisons made by Quickselect to find the smallest of a list of $n$ distinct numbers, and let $W_{n,1}$ be given by \eqref{eq:intro.def.wnm}. Then for all $n \ge 1$
	\beas 
	d_1(W_{n,1},D) \le \frac{8\log (n/2)+10}{n}.
	\enas
\end{theorem}

\begin{theorem} \label{thm:Quickselect.m}
	Let $m \ge 2$ and $C_{n,m}$ the number of comparisons made by Quickselect to find the $m^{th}$ smallest element of a list of $n$ distinct numbers, and let $W_{n,m}$ be given by \eqref{eq:intro.def.wnm}. Then for all $n \ge m$
	\beas
	d_1(W_{n,m},D) \le \frac{(46m+8) \log (n/m) + 54m+8}{n}.
	\enas
\end{theorem}

That the bounds in Theorems  \ref{thm:Quickselect} and \ref{thm:Quickselect.m} are tight in the $\log n/n$ order for $m \not =2$ is a consequence of the following result; in the following, we let $h_n=\sum_{1 \le k \le n}1/k$ for $n \ge 1$.

\begin{theorem} \label{thm:lower.bound}
	For all $m \ge 1$,
	\beas
	d_1(W_{n,m},D) \ge \frac{2(|m-2| \log n - |(m+2)h_m-3|)}{n} \qm{for all $n \ge m$.}
	\enas
\end{theorem}

We note that in the case $m=1$ the lower bound simplifies to $2 \log n/n$. That our method, where we focus only on the expectation $E[C_{n,m}]$ to achieve our lower bound, does not succeed in the case  $m=2$ is explained by the lack of the term $h_n$ on the right hand side of \eqref{eq:2n-4+2/n.m=2}. Theorem \ref{thm:lower.bound} is shown using the following exact expression for the expected running time of Quickselect; see also Section 6 of \cite{MR2672478}.

\begin{theorem}[Knuth \cite{MR0403310}] \label{thm:ECnm}
	Let $C_{n,m}$ be the number of comparisons made by Quickselect to locate the $m^{th}$ smallest of $n$ distinct numbers. Then for all  $n \ge m \ge 1$
	\bea \label{eq:EC.Theorem}
	E[C_{n,m}]= 2[n + 3 + (n + 1)h_n-(m + 2)h_m-(n -m + 3)h_{n-m+1}]. 
	\ena
\end{theorem}

In particular,
\bea \label{eq:2n-2h_n}
E[C_{n,1}]&=&2n-2h_n\\
E[C_{n,2}]&=&2n-4+\frac{2}{n}\label{eq:2n-4+2/n.m=2} \\
E[C_{n,3}]&=&2n-\frac{25}{3}+2h_n+\frac{2}{n-1}\qm{and}\nonumber\\\nonumber
E[C_{n,4}]&=&2n-13+4h_n-\frac{2}{n}+\frac{2}{n-2}.
\ena

Theorems \ref{thm:Quickselect} and \ref{thm:Quickselect.m} are derived by applying Theorem \ref{thm:DirectCoupling} that quantifies the if direction of the fixed point property \eqref{eq:D.is.U(D+1)} in the 
Wasserstein, or $d_1$ metric between two random variables $X$ and $Y$, given by 
\bea \label{def:Wass.is.sup}
d_1(X,Y) = \sup_{h \in {\rm Lip}_1}|Eh(X)-Eh(Y)| 
\qmq{where} {\rm Lip}_1 = \{h: |h(y)-h(x)|\le |y-x|\}.
\ena
On the left hand side of \eqref{def:Wass.is.sup} we have chosen to write $d_1(X,Y)$, rather than the technically correct expression $d_1({\cal L}(X),{\cal L}(Y))$, only for notational convenience.

\begin{theorem} \label{thm:DirectCoupling}
	Let $W$ be a non-negative random variable with finite mean, let $\theta>0$, and let the law of $W^*$ be given by \eqref{eq:W*transform}. Then 
	\beas 
	d_1(W,D_\theta) \le (1+\theta) d_1(W^*,W).
	\enas
\end{theorem}

As the Wasserstein distance also satisfies 
\bea \label{eq:d1.as.inf.over.couplings}
d_1(X,Y) = \inf E|X-Y|
\ena
where the infimum is over all couplings $(X,Y)$ having the given marginals, and is achieved here (see \cite{MR1105086}, for instance), Theorem \ref{thm:DirectCoupling} implies that
\bea \label{eq:d1bound.of.Theorem.coup}
d_1(W,D_\theta) \le (1+\theta)E|W^*-W|
\ena
for any non-negative random variable $W$ with finite mean, and $W^*$ defined on a common space having the ${\cal D}_\theta$-bias distribution of $W$.

In Section \ref{sec:quick} we detail the workings of the Quickselect algorithm and prove Theorems \ref{thm:Quickselect} and \ref{thm:Quickselect.m} by applying Theorem \ref{thm:DirectCoupling}, which is proved in Section \ref{sec:DirectCoupling}. The proof of Theorem \ref{thm:lower.bound} appears in Section
\ref{sec:proof.exact.ECnm}.

In related work, \cite{MR1932675} considers the Quicksort method, which produces a fully sorted list, and \cite{MR3248352} obtains distributional bounds for the running time of a variation of Quickselect to a non-Dickman approximand; compare its characterizion in (1.4) there to \eqref{eq:W*transform} here.

\section{The Quickselect Method and the Proofs of Theorems \ref{thm:Quickselect} and \ref{thm:Quickselect.m}} \label{sec:quick}
In this section we apply Theorem \ref{thm:DirectCoupling} to obtain the bounds in Theorems \ref{thm:Quickselect} and \ref{thm:Quickselect.m} on the error of the Dickman approximation for the distribution of $W_{n,m}$ in \eqref{eq:intro.def.wnm}, the properly normalized running time of the Quickselect algorithm for finding the $m^{th}$ smallest element of a list of $n$ distinct numbers. When the value of $m$ is clear from context, we will write $C_n$ for $C_{n,m}$.

\subsection{Quickselect: the case $m=1$}
In this section we prove Theorem \ref{thm:Quickselect} for the distribution of the number $C_n$ of comparisons that Quickselect requires to locate the smallest element of a list of $n$ distinct numbers.  Clearly, a list of size zero requires no comparisons, hence $C_0=0$. For $n \ge 1$, the procedure requires the $n-1$ comparisons of the pivot to every other element at the first stage, followed by the cost of processing the left subtree, which may be empty. Since the pivot is chosen uniformly, we obtain the stochastic recursion
\bea \label{eq:Cn.recursion}
C_n=n-1+C_{V_1} \qmq{for $n \ge 1$, with boundary condition $C_0=0$,}
\ena
where $V_1$, the size of the left subtree, is a discrete uniform variable on $\{0,\ldots,n-1\}$. From \eqref{eq:Cn.recursion} we see that $C_1=0$ and $C_2=1$ a.s., and that non-trivial distributions arise for $n \ge 3$.

Before proceeding to the proof of the theorem we describe how for all $n \ge 1$ we may write $C_n$ as a function $C(n;\textbf{U}_1)$ with
\beas
\textbf{U}_k=(U_k,U_{k+1},\ldots) \qmq{for $k \ge 1$,}
\enas
and $U_1,U_2,\ldots$ a sequence of i.i.d.\  uniform variables on $[0,1]$. 
Consider the initial list of size $V_0=n$ as making up the left subtree at stage 0. At stage $k \ge 1$, given a non-null left subtree from the previous stage of size $V_{k-1}$, a new left subtree of size
\bea \label{def:Uk}
V_k=\lfloor V_{k-1} U_k \rfloor \qm{for $k \ge 1$}
\ena
results by choosing a pivot uniformly from the current left subtree. In particular, the conditional distribution of $V_k$ given $V_{k-1}$ satisfies $V_k \sim {\cal U}\{0,\ldots,V_{k-1}-1\}$. Rewriting \eqref{eq:Cn.recursion} in this notation we have
\bea \label{eq:Cn.V.recursion}
C(n;\textbf{U}_1)=n-1+C(\lfloor nU_1\rfloor ;\textbf{U}_2) \qmq{for $n \ge 1$, with} C(0;\textbf{U}_k)=0 \qmq{for all $k \ge 1$.}
\ena
As the size of each non-null left subtree decrements by at least one at each iteration, the value of $C_n$ will only depend on an initial subsequence of $\textbf{U}_1$ of length at most $n$.

We pause to prove a lemma that is needed in this and the following section.
\begin{lemma} \label{lem:clog.over.m}
	If for $c$ a non-negative number and $q$ a positive integer
	\bea \label{eq:en.le.c+ave}
	e_n \le c +\frac{1}{n}\sum_{u=q}^{n-1} e_u \qm{for all $n \ge q$,}
	\ena
	then
	\bea \label{eq:en.le.clog}
	e_n \le  c \log (en/q) \qmq{for $n\ge q$.}
	\ena
\end{lemma} 

\noindent \emph{Proof:} As \eqref{eq:en.le.c+ave} holds for $n=q$ we see that $e_q \le c$, verifying that the inequality in \eqref{eq:en.le.clog} holds at $q$. 
Assuming inequality \eqref{eq:en.le.c+ave} holds for $q \le u \le n-1$ for some $n \ge q+1$ we have
\begin{multline*}
e_n \le c +\frac{1}{n}\sum_{u=q}^{n-1} e_u \le c +\frac{c}{n}\sum_{u=q}^{n-1} \log(eu/q)
\le c+ \frac{c}{n} \int_q^n \log (eu/q) du \\
= c\left(1+ \frac{1}{n}\left[ u \log (eu/q)-u\right]\Bvert_q^n \right)
= c \left( 1+\frac{1}{n}\left[ n \log (en/q)-n\right]\right)=c \log (en/q),
\end{multline*}
completing the inductive step, and the proof. \bbox \\[1ex]

We now prove Theorem \ref{thm:Quickselect}. In the proof, we use Lemmas \ref{lem:offby1} and \ref{lem:floorfloor}, which appear with their proofs at end of this section. \\

\noindent \emph{Proof of Theorem \ref{thm:Quickselect}:} 
Take $n \ge 1$. With $V_k$ as in \eqref{def:Uk}, by \eqref{eq:Cn.V.recursion} the variable $W_n$ as given by \eqref{eq:intro.def.wnm} satisfies 
\beas 
W_n=\frac{1}{n}C(n;\textbf{U}_1)-1 = \frac{1}{n}(n-1+C(V_1;\textbf{U}_2))-1
= \frac{1}{n}(C(V_1;\textbf{U}_2)-1).
\enas

We now construct a variable with the $W_n^*$ distribution by first constructing $W_n'$ having the $W_n$ distribution. As $\textbf{U}_1$ and $\textbf{U}_2$ are equidistributed,
\beas
W_n':=\frac{1}{n}C(n,\textbf{U}_2)-1 =_d \frac{1}{n}C(n,\textbf{U}_1)-1=W_n,
\enas
and hence 
\beas 
W_n^*:= U_1(W_n'+1)= \frac{1}{n}U_1C(n;\textbf{U}_2)
\enas
has the ${\cal D}$-bias distribution by \eqref{eq:W*transform}. The difference
\beas
W_n^*-W_n=\frac{1}{n}\left(U_1C(n;\textbf{U}_2)-C(V_1;\textbf{U}_2)+1\right)
\enas
satisfies
\beas
nE|W_n^*-W_n|  \le e_n + 1, \qmq{where we set} e_k=E|U_1C(k;\textbf{U}_2)-C(\lfloor k U_1 \rfloor;\textbf{U}_2)|, \quad k \ge 0,
\enas
hence consequence \eqref{eq:d1bound.of.Theorem.coup} of Theorem \ref{thm:DirectCoupling} with $\theta=1$ yields
\bea \label{eq:d1.2E.m=1.proof}
d_1(W_n,D) \le 2 E|W_n^*-W_n| \le \frac{2}{n}(e_n+1). 
\ena

We claim that
\begin{multline*}
e_n = E|U_1C(n;\textbf{U}_2)- C(\lfloor nU_1 \rfloor;\textbf{U}_2)|\\
\le  E|U_1(n-1)-
\lfloor nU_1 \rfloor+1|+ E|U_1C(\lfloor n U_2 \rfloor;\textbf{U}_3)-C(\lfloor \lfloor nU_1 \rfloor U_2 \rfloor;\textbf{U}_3)|.
\end{multline*}
When $\lfloor nU_1 \rfloor \ge 1$ this inequality follows from using the basic recursion \eqref{eq:Cn.V.recursion} on both terms forming the difference that defines $e_n$, followed by applying the triangle inequality, and is easily verified to hold directly in the case $\lfloor nU_1 \rfloor =0$ by applying \eqref{eq:Cn.V.recursion} only on the first term of that difference, noting the second one in this case is zero. Now using that $|u(n-1)-\lfloor nu \rfloor + 1| \le 2$ for all $u \in [0,1]$,  we obtain
\begin{multline} \label{eq:en.le.2.etc}
e_n \le 2 + E|U_1C(\lfloor nU_2 \rfloor;\textbf{U}_3)-  C(\lfloor \lfloor nU_1 \rfloor U_2\rfloor;\textbf{U}_3)|\\
\le 2+ E|U_1C(\lfloor nU_2 \rfloor;\textbf{U}_3)-  C(\lfloor \lfloor nU_2 \rfloor U_1 \rfloor;\textbf{U}_3)|+E|C(\lfloor \lfloor nU_2 \rfloor U_1 \rfloor;\textbf{U}_3)-C(\lfloor \lfloor nU_1 \rfloor U_2 \rfloor;\textbf{U}_3)|\\
= 2+ Ee_{\lfloor n U_2\rfloor}+E|C(\lfloor \lfloor nU_2 \rfloor U_1 \rfloor;\textbf{U}_3)-C(\lfloor\lfloor nU_1 \rfloor  U_2 \rfloor;\textbf{U}_3)|\\
\le 4 + Ee_{\lfloor n U_2\rfloor}.
\end{multline}
For the final term, the inequality
\beas
E|C(\lfloor \lfloor nU_2 \rfloor U_1 \rfloor;\textbf{U}_3)-C(\lfloor\lfloor nU_1 \rfloor  U_2 \rfloor;\textbf{U}_3)| \le 2 \qmq{for all $n \ge 0$}
\enas  
follows by applying Lemma \ref{lem:offby1}, below, that shows that $| \lfloor U_1 \lfloor n U_2 \rfloor \rfloor - \lfloor U_2 \lfloor nU_1 \rfloor \rfloor | \le 1$ a.s, and Lemma \ref{lem:floorfloor}, also below, that shows that $E|C(p,\textbf{U}_3)-C(p-1,\textbf{U}_3)| \le 2$ for all $p \ge 1$.

Expanding the expectation in $Ee_{\lfloor n U_2 \rfloor}$ in \eqref{eq:en.le.2.etc}, using the fact that $\lfloor nU_2 \rfloor$ is uniformly distributed over $\{0,\ldots,n-1\}$ and that $e_0=e_1=0$ by virtue of $C_0=C_1=0$, we obtain
\beas 
e_n \le 4+\frac{1}{n}\sum_{u=0}^{n-1} e_u \le 4+\frac{1}{n}\sum_{u=2}^{n-1} e_u\qmq{for $n \ge 2$.}
\enas
As $e_1=0$ inequality \eqref{eq:d1.2E.m=1.proof} shows that the claim of the theorem holds for $n=1$. Applying Lemma \ref{lem:clog.over.m} with $c=4$ and $q=2$ shows that $e_n \le 4 \log (en/2)$ for $n \ge 2$, and substituting this bound into 
\eqref{eq:d1.2E.m=1.proof} and simplifying now completes the proof. \bbox

We now prove Lemmas \ref{lem:offby1} and \ref{lem:floorfloor}. 
\begin{lemma} \label{lem:offby1}
	For all $(u_1,u_2) \in [0,1)^2$ and $n \ge 0$, 
	\beas
	| \lfloor u_1 \lfloor n u_2 \rfloor \rfloor - \lfloor u_2 \lfloor nu_1 \rfloor \rfloor | \le 1.
	\enas
\end{lemma}

\noindent \emph{Proof:} 
Consider the case $n \ge 1$, as otherwise the claim is trivial. Let $s=\lfloor n u_1  \rfloor$ and $t=\lfloor n  u_2  \rfloor$, so that $(s,t) \in \{0,1,\ldots,n-1\}^2$ and
\beas
s \le nu_1 < (s+1) \qmq{and} t \le nu_2 < (t+1).
\enas

Then
\beas
\frac{st}{n} \le u_2 \lfloor nu_1 \rfloor < \frac{s(t+1)}{n} \qmq{and} 
\frac{st}{n} \le u_1 \lfloor nu_2 \rfloor < \frac{(s+1)t}{n}.
\enas
Taking the difference, 
\beas
| u_1 \lfloor nu_2 \rfloor  - u_2 \lfloor nu_1 \rfloor| < 
\frac{1}{n}\max\{s,t\}<1.
\enas
As the difference between $u_1 \lfloor nu_2 \rfloor$  and $u_2 \lfloor nu_1 \rfloor$ is less than 1, their integer parts can differ by at most 1. 
\bbox

To prove Lemma \ref{lem:floorfloor},  we will use the easily verified fact that 
\bea \label{eq:ps.and.ks}
0 \le \frac{k-1}{p-1} <\frac{k}{p}<\frac{k}{p-1} \le 1\qm{for $p \ge 2$ and $1 \le k \le p-1$,}
\ena
and for $u \in [0,1]$ that
\bea \label{eq:rnd.un.unm1}
(\lfloor (p-1)u \rfloor, \lfloor pu \rfloor)= 
\left\{
\begin{array}{cc}
	(k-1,k-1) & u \in  \left[ \frac{k-1}{p-1},\frac{k}{p}\right) \\
	(k-1,k) & u \in \left[ \frac{k}{p},\frac{k}{p-1}\right).
\end{array}
\right.
\ena
We will also require the following inequality that can be shown directly using induction.
\begin{lemma} \label{lem:2n/n+1}If $c \ge 0, f_1 =0$ and
	\beas 
	f_p \le c + \frac{1}{p(p-1)}\sum_{k=1}^{p-1}kf_k \qmq{for all $p \ge 2$}
	\enas	
	then $f_p \le 2c$ for all $p \ge 1$. 
\end{lemma}

\begin{lemma} \label{lem:floorfloor}
	For all $p \ge 1$
	\beas
	f_p:=E|C(p,\textbf{U}_1)-C(p-1,\textbf{U}_1)| \le 2. 
	\enas
\end{lemma}

\noindent \textbf{Proof:}
As $f_1=0$ we need only consider $p \ge 2$. In view of \eqref{eq:ps.and.ks} we may write
\begin{multline*}
f_p=E|C(p,\textbf{U}_1)-C(p-1,\textbf{U}_1)|\\
= \sum_{k=1}^{p-1} E\left[ |C(p,\textbf{U}_1)-C(p-1,\textbf{U}_1)| \Bvert U_1 \in \left[ \frac{k-1}{p-1},\frac{k}{p}\right) \right] P\left(U_1 \in \left[ \frac{k-1}{p-1},\frac{k}{p}\right)\right)\\
+ \sum_{k=1}^{p-1} E\left[ |C(p,\textbf{U}_1)-C(p-1,\textbf{U}_1)| \Bvert U_1 \in \left[ \frac{k}{p},\frac{k}{p-1}\right) \right] P\left(U_1 \in \left[ \frac{k}{p},\frac{k}{p-1}\right)\right).
\end{multline*}

We claim that the conditional expectation in the first sum is 1. Indeed, for the given range of $U_1$ the first case of \eqref{eq:rnd.un.unm1} yields $(\lfloor (p-1)U_1 \rfloor,\lfloor pU_1\rfloor)=(k-1,k-1)$, and now 
\eqref{eq:Cn.V.recursion} implies that on this event
\beas
C(p,\textbf{U}_1)-C(p-1,\textbf{U}_1)=p-1 + C(k-1,\textbf{U}_2)-(p-2+C(k-1,\textbf{U}_2))=1.
\enas
For the second sum, the second case of \eqref{eq:rnd.un.unm1} yields $(\lfloor (p-1)U_1 \rfloor,\lfloor pU_1\rfloor)=(k-1,k)$, and
\begin{multline*}
C(p,\textbf{U}_1)-C(p-1,\textbf{U}_1)=p-1 + C(k,\textbf{U}_2)-(p-2+C(k-1,\textbf{U}_2))\\
=1+C(k,\textbf{U}_2)-C(k-1,\textbf{U}_2).
\end{multline*}
Hence, 
\begin{multline*}
f_p = \sum_{k=1}^{p-1}P\left(U_1 \in \left[ \frac{k-1}{p-1},\frac{k}{p}\right)\right)\\
+\sum_{k=1}^{p-1}E\left[ |1+C(k,\textbf{U}_2)-C(k-1,\textbf{U}_2)| \Bvert U_1 \in \left[ \frac{k}{p},\frac{k}{p-1}\right)\right]P\left(U_1 \in \left[ \frac{k}{p},\frac{k}{p-1}\right)\right)\\
\le \sum_{k=1}^{p-1}P\left(U_1 \in \left[ \frac{k-1}{p-1},\frac{k}{p}\right)\right)+\sum_{k=1}^{p-1}\left(1+f_k\right)P\left(U_1 \in \left[ \frac{k}{p},\frac{k}{p-1}\right)\right)\\
= 1 +  \sum_{k=1}^{p-1}f_k P\left(U_1 \in \left[ \frac{k}{p},\frac{k}{p-1}\right)\right)= 1 +  \frac{1}{p(p-1)}\sum_{k=1}^{p-1}kf_k. 
\end{multline*}
Invoking Lemma \ref{lem:2n/n+1} with $c=1$ now completes the proof. \bbox

\subsection{Case of $m \ge 2$} 
In this section we prove Theorem \ref{thm:Quickselect.m} for the approximation of the distribution of the properly scaled value of the number $C_{n,m}$ of comparisons made by the Quickselect algorithm $Q_m$ to determine the $m^{th}$ smallest element of a list of $n$ distinct numbers in the case $m \ge 2$.

As the $m^{th}$ smallest element of the list does not exist when $n < m$, no comparisons are required and we may set $C_{n,m}=0$ over this range. 
In the non-trivial case $n \ge m$, $Q_m$ begins as for $m=1$ at the first stage by selecting a uniformly chosen pivot, giving rise, through $n-1$ comparisons to the pivot, to a left subtree of size $V_1$, uniformly distributed over $\{0,\ldots,n-1\}$, and a right subtree of size $n-1-V_1$. If $V_1 \ge m$ then the $m^{th}$ smallest element of the original list lies in the left subtree, and we may locate it by applying $Q_m$ to it. If $V_1=m-1$ then the pivot is the $m^{th}$ smallest element and the process stops. Otherwise $V_1<m-1$, and the $m^{th}$ smallest element is the $m-V_1-1^{st}$ smallest element in the right subtree, which we then locate by applying $Q_{m-V_1-1}$ to it. Hence, we obtain
\begin{multline} \label{eq:Cnm.rec}
C_{n,m}=0 \qmq{for $0 \le n \le m-1$, and} \\
C_{n,m}=n-1 +C_{V_1,m}\textbf{1}(V_1 \ge m) + C_{n-V_1-1,m-V_1-1}\textbf{1}(V_1<m-1)  \qm{for $n \ge m$.}
\end{multline}

We now develop a simple bound on the expectation $E[C_{n,m}]$. 
\begin{lemma} \label{lem:4n}
	Let $C_{n,m}$ be the number of Quickselect comparisons for locating the $m^{th}$ smallest element of a list of $n$ distinct numbers. Then for all $m \ge 1$, 
	\beas 
	E[C_{n,m}] \le 4n \qmq{for all $n \ge 0$.}
	\enas
\end{lemma}

\noindent \emph{Proof:} Recall $h_n$ is the harmonic series $\sum_{1 \le k \le n}1/k$ for $n \ge 1$. The claim is trivial unless $n \ge m$, and is also easily seen to be true for $m=1$ and $m=2$ using \eqref{eq:2n-2h_n} and \eqref{eq:2n-4+2/n.m=2}. Hence, we take $n \ge m \ge 3$.

For such $n$ and $m$, writing the difference between the two harmonic series below as a sum and separating out the last term for $j=m-2$, we have
\bea \label{eq:weighted.h.sum}
(n-m+3)(h_n-h_{n-m+1})= \sum_{j=0}^{m-3} \frac{n-m+3}{n-j}+1+\frac{1}{n-m+2} \le m,
\ena
the inequality holding since each ratio is bounded by 1. Hence, using the expression given for $E[C_{n,m}]$ in Theorem \ref{eq:EC.Theorem} and 
applying \eqref{eq:weighted.h.sum} to yield the first inequality below, we obtain the upper bound
\begin{multline*}
E[C_{n,m}]= 2[n + 3 + (n+1)h_n-(m + 2)h_m-(n-m+3)h_{n-m+1}]  \\
=  2[n + 3 + (m-2)h_n-(m + 2)h_m+(n-m+3)(h_n-h_{n-m+1})]\\
\le 2[n + 3 + (m-2)h_n-(m + 2)h_m+m]\\
=  2[n + 3 + (m+1)(h_n-h_m)-3h_n+m-h_m] \\
\le  2[n + (m+1)\left(\frac{n-m}{m+1}\right)-3(h_n-1)+m-h_m] \\
= 2[2n-3(h_n-1)-h_m] \le 4n. 
\end{multline*}
\bbox

Note that the indicator on the first term on the right hand side of \eqref{eq:Cnm.rec} may be dropped, due to the boundary condition there, on the line above. Now letting $C_m(n;\textbf{U}_1)$ be defined by rewriting \eqref{eq:Cnm.rec} as \eqref{eq:Cn.V.recursion} was derived from \eqref{eq:Cn.recursion},
we obtain
\begin{multline} \label{eq:Cn.V.recursion.m}
C_m(n;\textbf{U}_1)=0 \qmq{for $0 \le n \le m-1$, and otherwise} \\
C_m(n;\textbf{U}_1) 
=n-1  + C_m(\lfloor nU_1 \rfloor; \textbf{U}_2) \\+ C_{m-1-\lfloor nU_1 \rfloor}(n-1-\lfloor nU_1 \rfloor; \textbf{U}_2)\textbf{1}(\lfloor nU_1 \rfloor < m-1).
\end{multline}
We next provide the following result that parallels Lemma \ref{lem:floorfloor} for the case $m=1$.

\begin{lemma} \label{lem:Cm.diff.10+16m}
	For all $m \ge 2$ and $p \ge 1$
	\beas
	f_p:=E|C_m(p;\textbf{U}_1)-C_m(p-1;\textbf{U}_1)| \le 2+16m.
	\enas
\end{lemma}
\noindent \emph{Proof:} 
As $C_m(p;\textbf{U}_1)=0$ for all $0 \le p \le m-1$ we may take $p \ge m$. 
By the basic recursion \eqref{eq:Cn.V.recursion.m} we have
\begin{multline*}
C_m(p;\textbf{U}_1)-C_m(p-1;\textbf{U}_1)
= 
1+C_m(\lfloor pU_1 \rfloor;\textbf{U}_2)-C_m(\lfloor (p-1)U_1 \rfloor;\textbf{U}_2)\\
+ C_{m-1-\lfloor pU_1 \rfloor}(p-1-\lfloor pU_1 \rfloor; \textbf{U}_2)\textbf{1}(\lfloor pU_1 \rfloor < m-1)\\
-C_{m-1-\lfloor (p-1)U_1 \rfloor}(p-2-\lfloor (p-1)U_1 \rfloor; \textbf{U}_2)\textbf{1}(\lfloor (p-1)U_1 \rfloor < m-1)\\ 
:= 1
+(C_m(\lfloor pU_1 \rfloor;\textbf{U}_2)-C_m(\lfloor (p-1)U_1 \rfloor;\textbf{U}_2))+R.
\end{multline*}
Applying the triangle inequality and taking expectation yields
\bea \label{eq:fn.bound.linear.in.m}
f_p \le 1+E|C_m(\lfloor pU_1 \rfloor;\textbf{U}_2)-C_m(\lfloor (p-1)U_1 \rfloor;\textbf{U}_2)|+E|R|.
\ena
For the first expectation in \eqref{eq:fn.bound.linear.in.m}, by \eqref{eq:rnd.un.unm1} we have 
\begin{multline*}
E|C_m(\lfloor pU_1 \rfloor;\textbf{U}_2)-C_m(\lfloor (p-1)U_1 \rfloor;\textbf{U}_2)|\\
= \sum_{k=1}^{p-1}E\left[ |C_m(k-1,\textbf{U}_2)-C_m(k-1,\textbf{U}_2)|\Bvert U_1 \in \left[ \frac{k-1}{p-1},\frac{k}{p}\right) \right]P\left(U_1 \in \left[ \frac{k-1}{p-1},\frac{k}{p}\right)\right)\\
+\sum_{k=1}^{p-1}E\left[ |C_m(k,\textbf{U}_2)-C_m(k-1,\textbf{U}_2)|\Bvert U_1 \in \left[ \frac{k}{p},\frac{k}{p-1}\right) \right]P\left(U_1 \in \left[ \frac{k}{p},\frac{k}{p-1}\right)\right)\\
= \sum_{k=1}^{p-1} f_k P\left(U_1 \in \left[\frac{k}{p},\frac{k}{p-1}\right) \right)
= \frac{1}{p(p-1)}\sum_{k=1}^{p-1} kf_k.
\end{multline*}

Now applying Lemma \ref{lem:4n} on the first term of the remainder $R$, and using that $\lfloor pU_1 \rfloor \sim {\cal U}\{0,\ldots,p-1\}$, yields
\begin{multline*}
E[C_{m-1-\lfloor pU_1 \rfloor}(p-1-\lfloor pU_1 \rfloor; \textbf{U}_2)\textbf{1}(\lfloor pU_1 \rfloor < m-1)] \le \frac{4}{p}\sum_{k=0}^{m-2}(p-1-k) \\
\le  \frac{4}{p}(p-1)(m-1) \le 4m,
\end{multline*}
and replacing $p$ by $p-1$ we see that the same bound holds for the expectation of the final term of $R$. 

Substituting the bounds achieved into \eqref{eq:fn.bound.linear.in.m} we obtain
\bea \label{eq:m.case.invoke.lemma2n/n+1}
f_p \le 1+8m+\frac{1}{p(p-1)}\sum_{k=1}^{p-1} kf_k \qm{for all $p \ge m$.}
\ena
As $f_p=0$ for $1 \le p \le m-1$ inequality \eqref{eq:m.case.invoke.lemma2n/n+1} holds for all $p \ge 2$, and the conditions for invoking Lemma \ref{lem:2n/n+1} with
$c=1+8m$ are satisfied, yielding the desired conclusion. \bbox \\[1ex]

\noindent {\em Proof of Theorem \ref{thm:Quickselect.m}:} 
Let $n \ge m$. From \eqref{eq:intro.def.wnm} and \eqref{eq:Cn.V.recursion.m}, letting $V_1=\lfloor nU_1 \rfloor$, 
\begin{multline} \label{Dn:one.level.down.m}
W_n=\frac{1}{n}C_m(n;\textbf{U}_1)-1 \\= \frac{1}{n}(n-1+C_m(V_1;\textbf{U}_2)+C_{m-1-V_1}(n-1-V_1;\textbf{U}_2)\textbf{1}(V_1 < m-1))-1\\
= \frac{1}{n}(C_m(V_1;\textbf{U}_2)+C_{m-1-V_1}(n-1-V_1;\textbf{U}_2)\textbf{1}(V_1 < m-1)-1).
\end{multline}

We now construct a variable with the $W_n^*$ distribution. As $\textbf{U}_1$ and $\textbf{U}_2$ are equidistributed, $W_n'$ given by the first equality in \eqref{Dn:one.level.down.m} when substituting $\textbf{U}_2$ in place of $\textbf{U}_1$ has law ${\cal L}(W_n)$. Hence, by \eqref{eq:W*transform} with $\theta=1$, letting
\bea \label{eq:Dn*.quickselect}
W_n^*= U_1(W_n'+1)= \frac{1}{n}U_1C_m(n;\textbf{U}_2),
\ena
the pair $(W_n,W_n^*)$ is a coupling of a variable with the $W_n$ distribution to one with its Dickman ${\cal D}$-bias distribution. Applying consequence \eqref{eq:d1bound.of.Theorem.coup} of Theorem \ref{thm:DirectCoupling}, we obtain
\bea \label{eq:factor.of.2.theorem.m}
d_1(W_n,D) \le \frac{2}{n}f_n \qmq{where} f_n=nE|W_n^*-W_n|.
\ena
Letting
\beas 
e_n=E|U_1C_m(n;\textbf{U}_2)-C_m(\lfloor n U_1 \rfloor;\textbf{U}_2)|,
\enas
in view of \eqref{Dn:one.level.down.m} and
\eqref{eq:Dn*.quickselect}, and applying Lemma \ref{lem:4n} to bound expectations of the form $E[C_{n,m}]$ and that $V_1 \sim {\cal U}\{0,1,\ldots,n-1\}$, we obtain 
\begin{multline}
f_n=nE|W_n^*-W_n|\\
=E|U_1C_m(n;\textbf{U}_2)-C_m(V_1;\textbf{U}_2)-C_{m-1-V_1}(n-1-V_1;\textbf{U}_2)\textbf{1}(V_1 < m-1)+1|\\
\le e_n + E|C_{m-1-V_1}(n-1-V_1;\textbf{U}_2)\textbf{1}(V_1 < m-1)|+1\\
\le e_n + \frac{4}{n}\sum_{k=0}^{m-2} (n-1-k) +1\\
\le  e_n + \frac{4}{n}(n-1)(m-1)+1
\le e_n + 4m. \label{eq:gen.m.fn.le.2n+4m}
\end{multline}

To control $e_n$, invoke the basic recursion \eqref{eq:Cn.V.recursion.m} to write
\begin{multline*}
U_1C_m(n;\textbf{U}_2)
=U_1(n-1)+U_1C_m(\lfloor nU_2 \rfloor;\textbf{U}_3)\\+U_1C_{m-1-\lfloor nU_2 \rfloor}(n-1-\lfloor nU_2 \rfloor;\textbf{U}_3)\textbf{1}(\lfloor nU_2 \rfloor < m-1)\\
= U_1(n-1)+U_1C_m(\lfloor nU_2 \rfloor;\textbf{U}_3)+R_1
\end{multline*}
where
\beas
R_1=U_1C_{m-1-\lfloor nU_2 \rfloor}(n-1-\lfloor nU_2 \rfloor;\textbf{U}_3)\textbf{1}(\lfloor nU_2 \rfloor < m-1),
\enas
and similarly, 
\beas
C_m(\lfloor n U_1 \rfloor;\textbf{U}_2)
= (\lfloor nU_1 \rfloor-1)\textbf{1}(\lfloor n U_1 \rfloor \ge m) +  C_m(\lfloor \lfloor nU_1 \rfloor U_2 \rfloor;\textbf{U}_3) + R_2\\
= (\lfloor nU_1 \rfloor-1) +  C_m(\lfloor \lfloor nU_1 \rfloor U_2 \rfloor;\textbf{U}_3) + R_2 +R_3
\enas
where
\begin{multline*}
R_2=C_{m - 1-\lfloor  \lfloor nU_1 \rfloor U_2 \rfloor}( \lfloor nU_1 \rfloor-1  - \lfloor  \lfloor nU_1 \rfloor U_2 \rfloor;\textbf{U}_3)\textbf{1}(\lfloor \lfloor nU_1 \rfloor U_2 \rfloor < m-1,\lfloor nU_1 \rfloor \ge m), 
\end{multline*}
and
\beas
R_3=-(\lfloor nU_1 \rfloor-1)\textbf{1}(\lfloor n U_1 \rfloor \le m-1). 
\enas
Taking the expectation of the absolute difference and using that $|u(n-1)-\lfloor nu \rfloor + 1| \le 2$ for all $u \in [0,1]$, we obtain
\begin{multline}
e_n=E|U_1C_m(n;\textbf{U}_2)- C_m(\lfloor n U_1 \rfloor;\textbf{U}_2)| \\
\le E|U_1(n-1)-(\lfloor nU_1 \rfloor-1)|
+E|U_1C_m(\lfloor nU_2 \rfloor;\textbf{U}_3) \\-  C_m(\lfloor \lfloor nU_1 \rfloor U_2 \rfloor;\textbf{U}_3)|+E|R_1|+E|R_2|+E|R_3|\\
\le 2+E|U_1C_m(\lfloor nU_2 \rfloor;\textbf{U}_3) -  C_m(\lfloor \lfloor nU_1 \rfloor U_2 \rfloor;\textbf{U}_3)|+E|R_1|+E|R_2|+E|R_3|\\
\le 2+E|U_1C_m(\lfloor nU_2 \rfloor;\textbf{U}_3) -C_m(\lfloor \lfloor nU_2 \rfloor U_1 \rfloor;\textbf{U}_3)| \\+ E|C_m(\lfloor \lfloor nU_2 \rfloor U_1 \rfloor;\textbf{U}_3)-C_m(\lfloor \lfloor nU_1 \rfloor U_2 \rfloor;\textbf{U}_3)|
+E|R_1|+E|R_2|+E|R_3|.\label{eq:m.ge.2.R1to3}
\end{multline}

Lemmas \ref{lem:offby1} and \ref{lem:Cm.diff.10+16m} yield 
\bea \label{eq:gen.m.R0}
E|C_m(\lfloor \lfloor nU_2 \rfloor U_1 \rfloor;\textbf{U}_3)-C_m(\lfloor \lfloor nU_1 \rfloor U_2 \rfloor;\textbf{U}_3)| \le 2+16m.
\ena

For the first remainder term $R_1$, by Lemma \ref{lem:4n}, we have 
\begin{multline} \label{eq:gen.m.R1}
E|R_1| = \frac{1}{2}E[C_{m-1-\lfloor nU_2 \rfloor}(n-1-\lfloor nU_2 \rfloor;\textbf{U}_3)\textbf{1}(\lfloor nU_2 \rfloor \le m-2)]\\
\le \frac{2}{n}\sum_{k=0}^{m-2}(n-1-k)=\frac{2}{n}(n-1)(m-1) \le 2m. 
\end{multline}

For $R_2$, we condition on the event $\lfloor n U_1 \rfloor =k$ for $1 \le k \le n-1$, then further on $\lfloor kU_2 \rfloor=j$ for $0 \le j \le k-1$. We note the presence of $\lfloor n U_1 \rfloor \ge m$ in the indicator restricts $k \ge m \ge 2$ in this second step, where the values of $j$ are all equally likely with probability $1/k$. Applying Lemma \ref{lem:4n} then yields
\begin{multline} \label{eq:gen.m.R2}
E|R_2| = E[C_{m-1  - \lfloor  \lfloor nU_1 \rfloor U_2 \rfloor}( \lfloor nU_1 \rfloor-1  - \lfloor  \lfloor nU_1 \rfloor U_2 \rfloor;\textbf{U}_3)\textbf{1}(\lfloor \lfloor nU_1 \rfloor U_2 \rfloor < m-1, \lfloor nU_1 \rfloor \ge m)]\\
\le \frac{4}{n}\sum_{k=m}^{n-1}\frac{1}{k}\sum_{j=0}^{m-2} (k-1-j) \le \frac{4m}{n}\sum_{k=1}^{n-1}\frac{1}{k}(k-1)\le \frac{4m}{n}(n-1) \le 4m. 
\end{multline}

As $R_3$ satisfies
\bea \label{eq:gen.m.R3}
E|R_3|  = E|(\lfloor nU_1 \rfloor-1)\textbf{1}(\lfloor n U_1 \rfloor \le m-1)| \le m,
\ena
substituting the bounds \eqref{eq:gen.m.R0}-\eqref{eq:gen.m.R3} into \eqref{eq:m.ge.2.R1to3} yields that, for all $n \ge m$, 
\begin{multline*}
e_n \le 4+23m +E|U_1C_m(\lfloor nU_2 \rfloor;\textbf{U}_3) -  C_m(\lfloor \lfloor nU_2 \rfloor U_1 \rfloor;\textbf{U}_3)|\\
=4+23m+\frac{1}{n}\sum_{k=0}^{n-1}e_k =4+23m+\frac{1}{n}\sum_{k=m}^{n-1}e_k,
\end{multline*}
where the final equality follows by noting that $C(k;\textbf{U}_1)=0$ for $k \le m-1$. Applying Lemma \ref{lem:clog.over.m} yields that, for all $n \ge m$, 
\beas
e_n \le (4 + 23m) \log (ne/m), 
\enas
and now from \eqref{eq:gen.m.fn.le.2n+4m} we conclude
\beas
f_n \le e_n + 4m = (4 + 23m) \log (ne/m)+4m.
\enas Substitution into \eqref{eq:factor.of.2.theorem.m}, and simplification, yields the claim. \bbox

\section{Proof of Theorem 1.5} \label{sec:DirectCoupling}
Theorem \ref{thm:DirectCoupling} was originally proven using Stein's method in \cite{goldstein2017non}, but \cite{RN} offered the following much simpler approach. \\

\noindent \emph{Proof:} Let $U \sim {\cal U}[0,1]$ be independent of the pair $(W,D_\theta)$, which are constructed on the same space so as to achieve the infimum in \eqref{eq:d1.as.inf.over.couplings}. Then, as $D_\theta=_d D_\theta^*$, 
\beas
d_1(W^*,D_\theta) = d_1(U^{1/\theta}(W+1),U^{1/\theta}(D_\theta+1)) 
\le E[U^{1/\theta}|W-D_\theta|] = \frac{\theta}{\theta+1}d_1(W,D_\theta).
\enas
Now, by the triangle inequality, 
\beas
d_1(W,D_\theta) \le d_1(W,W^*)+ d_1(W^*,D_\theta) \le d_1(W,W^*)+\frac{\theta}{\theta+1}d_1(W,D_\theta).
\enas
Rearranging the inequality yields the claimed bound. \bbox

\section{Proof of Theorem 1.3}
\label{sec:proof.exact.ECnm}

\noindent We now apply Theorem \ref{thm:ECnm} to prove Theorem \ref{thm:lower.bound}.\\

\noindent {\em Proof of Theorem \ref{thm:lower.bound}.} Since $f(x)=x$ is an element of ${\rm Lip}_1$, expression \eqref{def:Wass.is.sup} for 
the Wasserstein distance yields that
\beas
d_1(W_{n,m},D) \ge |E[W_{n,m}]-E[D]| =|E[W_{n,m}]-1|=\left|\frac{1}{n}E[C_{n,m}]-2\right|,
\enas
applying \eqref{eq:intro.def.wnm} and that (see e.g. \cite{MR1918722}) $E[D]=1$.

Now, slightly rewriting the equality in \eqref{eq:EC.Theorem} as 
\begin{align*}
E[C_{n,m}] = 2\left[n+3+(m-2) h_n - (m+2)h_m  +(n-m+3)(h_n-h_{n-m+1}) \right]
\end{align*}
for $m>2$ we have
\begin{align*}
\frac{1}{n}E[C_{n,m}]-2  \ge \frac{2[(m-2)h_n-(m+2)h_m+3]}{n}  \ge \frac{2[(m-2) \log n  -|(m+2)h_m-3|]}{n},
\end{align*}
using $h_n > \log n$. Hence, the claim of Theorem \ref{thm:lower.bound} holds for $m>2$. We see the claim of Theorem  also holds for $m=1$ by using the form \eqref{eq:2n-2h_n}, which yields $|E[C_{n,m}/n-2|=2h_n/n$, noting that in this case $(m+2)h_m-3=0$.  \bbox \\[1ex]

\noindent {\bf Acknowledgement} The author thanks Ralph Neininger for his vast simplification of the previous proof of Theorem \ref{thm:DirectCoupling} in the preprint \cite{goldstein2017non}, as well as for the suggestion for obtaining the lower bounds as achieved in Theorem \ref{thm:lower.bound}. The author also sincerely thanks two reviewers whose suggestions and observations were extremely valuable, which included pointing out that Theorem \ref{thm:ECnm} is a known result due to Knuth, and a simplification of Lemma \ref{lem:4n}.

\end{document}